\documentclass[11pt,draft,twoside]{article}
\usepackage[cp1250]{inputenc}
\usepackage{latexsym,amssymb,amsmath}
\title{The commutative Moufang loops with
maximum conditions for subloops}
\date{}
\author{A. Babiy, N. Sandu}
\begin{document}
\maketitle

\begin{abstract}

It is proved that the maximum condition for subloops in a
commutative Moufang loop $Q$ is equivalent with the conditions of
finite generating of different subloops of the loop $Q$ and
different subgroups of the multiplication group of the loop $Q$.
An analogue equivalence is set for the commutative Moufang
$ZA$-loops.

\ Classification: {\ 20N05}

\ Keywords and phrases: {\ commutative Moufang loop,
multiplication group of loop, maximum condition for subloops.}

\end{abstract}

It is said that maximum condition (respect. minimum condition) for
the subalgebras with the property $\alpha$ holds  in an algebra
$A$ if any ascending (respect. descending) system of subalgebras
with the property $\alpha$ $A_1 \subseteq A_2 \subseteq \ldots$
(respect. $A_1 \supseteq A_2 \supseteq \ldots$ ) break, i.e. $A_n
= A_{n+1} \ldots $ for a certain $n$. It is well known that the
fulfillment of the maximum condition for subalgebras of an
arbitrary algebra is equivalent to the fact that both the algebra
and any of its subalgebras are finitely generated.

Commutative Moufang loops (CML's) with maximum condition for
sub\-loops is considered in this paper. It is proved that for a
non-associati\-ve CML $Q$ this condition is equivalent to one of
the following equivalent conditions: a) if  $Q$ contains a
centrally nilpotent subloop of class $n$, then all its subloops of
this type are finitely generated; b) if  $Q$ contains a centrally
solvable subloop of class $s$, then all its subloops of this type
are finitely generated; c) all invariant subloops of  $Q$ are
finitely generated; d) all non-invariant associative subloops of
$Q$ are finitely generated; e) at least one maximal associative
subloop of $Q$ is finitely generated. This list is completed with
the condition of finite generating of various subgroups of the
multiplication group of  $Q$. If $Q$ is a $ZA$-loop, then the list
a) - e) is completed with the condition of finite generating of
the center of  $Q$, as well  with the condition of finite
generating of other subloops of  $Q$ and various subgroups of the
multiplication group of  $Q$.

It is worth  mentioning that the following statement is proved in
[1, 2].
\smallskip\\
\textbf{Lemma 1.} \textit{The following conditions are equivalent
for an arbitrary CML $Q$:}

\textit{1) $Q$ is finitely generated;}

\textit{2) the maximum condition for subloops holds in  $Q$.}

In [2] the list a) - e) is completed with equivalent statements:
h) the CML $Q$ satisfies the maximum condition for invariant
subloops; i) the CML $Q$ is a subdirect product of a finite CML of
exponent 3 and a finitely generated abelian group; j) the CML $Q$
possesses a finite central series, whose factors are cyclic groups
of simple or infinite order.

Let us bring some notions and results on the theory of commutative
Moufang loops, needed for the further research.

A \textit{commutative Moufang loop} (CML's) is characterized by
the identity $ x^2\cdot yz = xy\cdot xz$. The
\textit{multiplication group} $\frak M(Q)$ of a CML $Q$ is the
group generated by all the \textit{translations} $L(x)$, where
$L(x)y = xy$. The subgroup $I(Q)$ of the group $\frak M(Q)$,
generated by all the \textit{inner mappings} $L(x,y) =
L(xy)^{-1}L(x)L(y)$ is called the \textit{inner mapping group} of
the CLM $Q$. The subloop $H$ of a CML $Q$ is called
\textit{normal} (\textit{invariant}) in $Q$, if $I(Q)H = H$.
\smallskip\\
\textbf{Lemma 2} [3]. \textit{Let $Q$ be a commutative Moufang
loop with the multiplication group $\frak M$. Then $\frak
M/Z(\frak M)$, where $Z(\frak M)$ is the centre of the group
$\frak M$, and $\frak M^{\prime} = (\frak M,\frak M)$ are locally
finite $3$-groups and will be finite if $Q$ is finitely
generated}.
\smallskip\\
\textbf{Lemma 3}. \textit{The multiplication group $\frak M$ of an
arbitrary CML is locally nilpotent}.
\smallskip\\
\textbf{Proof.} Let $\overline{\frak N}$ be the image of finitely
generated subgroup of group $\frak M$  under the homomorphism
$\frak M \rightarrow \frak M/Z(\frak M)$. It follows from Lemma 2
that $\overline{\frak N}$ is a finite $3$-group, therefore it is
nilpotent. Let us write $\overline{\frak N}$ in the form $\frak N
Z(\frak M)/Z(\frak M)$. We have $\frak N Z(\frak M)/Z(\frak M)
\cong$ $\frak N/(\frak N \cap Z(\frak M))$. It is obvious that
$\frak N \cap Z(\frak M) \subseteq Z(\frak N)$. Then
$$ \frak N/Z(\frak N) \cong (\frak N/(\frak N \cap Z(\frak
M)))/(Z(\frak N)/(\frak N \cap Z(\frak M))). $$

Therefore $\frak N/Z(\frak N)$ is nilpotent, as a homomorphic
image of the nilpotent group $\frak N/(\frak N \cap Z(\frak M))$.
Then the group $\frak N$ is nilpotent as well. Consequently, the
group $\frak M$ is locally nilpotent, as required.

The \textit{center} $Z(Q)$ of a CML $Q$ is an invariant subloop
$Z(Q) = \{x \in Q \vert x\cdot yz = xy\cdot z \forall y,z \in
Q\}$.
\smallskip\\
\textbf{Lemma 4} [3]. \textit{Quotient loop $Q/Z(Q)$ of an
arbitrary CML $Q$ on its center $Z(Q)$ has the exponent three.}
\smallskip\\
\textbf{Lemma 5} [3]. \textit{A periodic CML is locally finite.}

The associator $(a,b,c)$ of the elements $a, b, c$ in CML $Q$ is
defined by the equality $ab\cdot c = (a\cdot bc)(a,b,c)$. We
denote by $Q_i$ (respect. $Q^{(i)}$) the subloop of the CML $Q$,
generated by all associators of the form
$(x_1,x_2,\ldots,x_{2i+1})$ (respect.
$(x_1,\ldots,x_{3^i})^{(i)}$), where
$(x_1,\ldots,x_{2i-1},x_{2i},x_{2i+1}) =
((x_1,\ldots,x_{2i-1}),x_{2i},x_{2i+1})$ (respect.
$(x_1,\ldots,x_{3^i})^{(i)} = ((x_1,\ldots,x_{3^{
i-1}})^{(i-1)},(x_{3^{i-1}+1},\ldots,x_{2\cdot 3^{i-1}})^{(i-1)},
\break (x_{2\cdot 3^{i-1}+1},\ldots, x_{3^i})^{(i-1)})$, where
$(x_1,x_2,x_3)^{(1)} = (x_1,x_2,x_3)$ . The series of normal
subloops $1 = Q_0 \subseteq Q_1 \subseteq \ldots \subseteq Q_i
\subseteq \ldots$ (respect. $1 = Q^{(o)} \subseteq Q^{(1)}
\subseteq \ldots \break \ldots \subseteq Q^{(i)} \subseteq \ldots
$) is called the \textit{lower central series} (respect.
\textit{derived series}) of the CML $Q$. We will also use for
associator loop the designation $Q^{(1)} = Q'$.

A CML $Q$ is \textit{centrally nilpotent} (respect.
\textit{centrally solvable}) \textit{of class $n$} if and only if
its lower central series (respect. derived series) has the form $1
\subset Q_1 \subset \ldots \subset Q_n = Q$ (respect. $1 \subset
Q^{(1)} \subset \ldots \subset Q^{(n)} = Q$) [3].

An \textit{ascending central series} of CML $Q$ is a linearly
ordered by the inclusion system

$$1 = Q_0 \subseteq Q_1 \subseteq \ldots \subseteq Q_{\alpha}
\subseteq \ldots \subseteq Q_{\gamma} = Q$$ of invariant subloops
of  $Q$, satisfying the conditions:

1) $Q_{\alpha} = \sum_{\beta < \alpha} Q_{\beta}$ for limit
ordinal $\alpha$;

2) $Q_{\alpha + 1}/Q_{\alpha} \subseteq Z(Q/Q_{\alpha})$.

A CML, possessing an ascending central series is called
\textit{$ZA$-loop}. If the ascending central series of CML is
finite, then it is centrally nilpotent [3].

We will often use the following statements in our further proofs.
\smallskip\\
\textbf{Lemma 6}. \textit{The following statements are equivalent
for an arbitrary CML $Q$:}

\textit{1)  $Q$ satisfies the minimum condition for subloops;}

\textit{2)  $Q$ is a direct product of a finite number of
quasicyclic groups, belonging to the center  CML $Q$, and a finite
CML;}

\textit{3)  $Q$ satisfies the minimum condition for invariant
subloops;}

\textit{4)  $Q$ satisfies the minimum condition for non-invariant
associative subloops;}

\textit{5) if  $Q$ contains a centrally nilpotent subloop of class
$n$, then it satisfies the minimum condition for centrally
nilpotent subloops of class $n$;}

\textit{6) if  $Q$ contains a centrally solvable subloop of class
$s$, then it satisfies the minimum condition for centrally
solvable subloops of class $s$;}

\textit{7) at least one maximal associative subloop of  $Q$
satisfies the minimum condition for subloops.}

The equivalence of conditions 1), 2), 3) is proved in [4], the
equivalence of conditions 1), 4), 5), 6) is proved in [5] and the
equivalence of conditions 1), 7) is proved in [6].
\smallskip\\
\textbf{Lemma 7} [4]. \textit{The following statements are
equivalent for an arbitrary non-associative CML $Q$ with a
multiplication group $\frak M$:}

\textit{1)  $Q$ satisfies the minimum condition for subloops;}

\textit{2)  $\frak M$ satisfies the minimum condition for
subgroup;}

\textit{3)  $\frak M$ is a product of a finite number of
quasicyclic groups, lying in the center of  $\frak M$, and a
finite group;}

\textit{4)  $\frak M$ satisfies the minimum condition for
invariant subgroup;}

\textit{5) at least one maximal abelian subgroup of  $\frak M$
satisfies the minimum condition for subgroups;}

\textit{6) if  $\frak M$ contains a nilpotent subgroup of  class
$n$, then $\frak M$ satisfies the minimum condition for nilpotent
subgroups of  class $n$.}

\textit{7) if  $\frak M$ contains a solvable subgroup of  class
$s$, then $\frak M$ satisfies the minimum condition for solvable
subgroups of  class $s$.}
\smallskip\\
\textbf{Lemma 8} [4]. \textit{If the center $Z(Q)$ of a
commutative Moufang $ZA$-loop $Q$ satisfies the minimum condition
for subloops, then  $Q$ satisfies the minimum condition for
subloop itself.}

Let us now consider an arbitrary non-periodic CML $Q$. Let $Q^3 =
\{x^3 \vert x \in Q\}$. CML is di-associative [3], then it is easy
to show that $Q^3$ is a subloop. It follows from Lemma 4 that $Q^3
\subseteq Z(Q)$, where $Z(Q)$ is the center of CML $Q$, therefore
$Q^3$ is an invariant subloop of  $Q$. Let us suppose that the
subloop $Z(Q)$ is finitely generated. Then the abelian group $Q^3$
is also finitely generated. Therefore it decomposes into a direct
product of cyclic groups $Q^3 = <r_1> \times \ldots \times <r_k>
\times <s_1> \times \ldots \times <s_m> = <R> \times <S>$, where
$<r_i>$ are cyclic groups of infinite order, $<s_j>$ are finite
cyclic groups [7]. T he group $R$ is free abelian, therefore it is
without torsion. It is shown in [3] that the associator loop $Q'$
has the exponent three, then

$$R \cap Q' = \{1\}. \eqno{(1)}$$
\smallskip\\
\textbf{Lemma 9.} \textit{Let $Q$ be a CML, $R$ be its subloop,
which is considered above, and let $\overline H$ be a subloop of
CML $Q/R = \overline Q$. The subloop $\overline H$ satisfies one
of the properties: 1)  $\overline H$ is centrally nilpotent of
class $n$; 2)  $\overline H$ is centrally solvable of class $s$;
3) $\overline H$ is a maximal associative subloop of CML
$\overline Q$; 4) $\overline H$ is the center of CML $\overline
Q$; 5) $\overline H$ is a non-invariant subloop of CML $Q$; 6)
$\overline H$ is an invariant subloop of CML $\overline Q$ if and
only if the inverse image $H$ of subloop $\overline H$ has the
same property as the subloop $\overline H$, under the homomorphism
$\varphi: Q \rightarrow Q/R$.}
\smallskip\\
\textbf{Proof.} Let us suppose that subloop $\overline H$ is
centrally nilpotent of class $n$. Let $h_1, h_2, \ldots, h_{2n+1}$
be arbitrary elements from $H$. Let us denote $\varphi(h_i) =
\overline h_i, \varphi (1) = \overline 1$. Then $\overline h_i =
h_iR, \overline 1 = R$. We have $(\overline h_1,\overline h_2,
\ldots, \overline h_{2n+1}) = \overline 1, (h_1R,h_2R, \ldots,
h_{2n+1}R) = R$. But $R \subseteq Z(Q)$. Therefore, if $u \in R$,
then $(au,b,c) = (a,b,c)$ for any elements $a, b, c \in Q$. Then
$(h_1,h_2, \ldots, h_{2n+1}) = r$, where $r \in R$. It follows
from (1) that $r = 1$. We have obtained that $(h_1,h_2, \ldots,
h_{2n+1}) = 1$, i.e. the subloop $H$ is centrally nilpotent of
class $n$.

Conversely, let us suppose that the subloop $H$ is centrally
nilpotent of class $n$. Then there exist such elements $h_1, h_2,
\ldots, h_{2n-1}$ from $H$ that $(h_1,h_2, \ldots, h_{2n-1}) \neq
1$. It follows from (1) that $(h_1,h_2, \ldots, h_{2n-1}) \notin
R$. Therefore $(\overline h_1,\overline h_2, \ldots, \overline
h_{2n-1}) \notin \overline 1$. Consequently $\overline H$, as
homomorphic image of subloop $H$, will be a centrally nilpotent
subloop of class $n$. It proves the statement 1). The statement 2)
is proved by analogy.

Let us now suppose that $\overline H$ is a maximal associative
subloop of CML $\overline Q$ and the inverse image $H$ is not a
maximal associative subloop of CML $Q$. Then there exists such an
element $a \notin H$, that $(a,h_1,h_2) = 1$ for all $h_1, h_2 \in
H$. Obviously $R \subseteq H$. Then $\varphi a = \overline a
\notin \overline H$ and $(\overline a,\overline h_1,\overline h_2)
= \overline 1$ for all $\overline h_1, \overline h_2 \in \overline
H$. We have obtained that the non-associative subloop $<\overline
a, \overline H>$, generated by the set $\{\overline a, \overline
H\}$, strictly contains $\overline H$, i.e. $\overline H$ is not a
maximal associative subloop of CML $\overline Q$. Contradiction.
Consequently, $H$ is a maximal associative subloop of CML $Q$.

Conversely, let us suppose that $H$ is a maximal associative
subloop of CML $Q$ and $\overline H$ is not a maximal associative
subloop of CML $\overline Q$. Then there exists such an element
$\overline a \notin \overline H$, that $(\overline a,\overline
h_1,\overline h_2) = \overline 1$ for all $\overline h_1,
\overline h_2 \in \overline H$. We have obtained that
$(aR,h_1R,h_2R) = R$ for all $h_1, h_2 \in H$. As $R \subseteq
Z(Q)$, then $(a,h_1,h_2) = r$, where $r \in R$. It follows from
(1) that $r = 1$, therefore $(a,h_1,h_2) = 1$ for all $h_1, h_2
\in H$ and $a \notin H$. It means that subloop $H$ is strictly
contained in the associative subloop $<a, H>$. We have obtained a
contradiction with the fact that subloop $H$ is a maximal
associative subloop. This proves statement 3). Statement 4) is
proved by analogy.

Statements 5), 6) follow from the fact that the natural
homomorphism $Q \rightarrow Q/R$ sets a one-to-one mapping between
all non-invariant (respect. invariant) subloops of CML $Q$, with
contained $R$, and all non-invariant (respect. invariant) subloops
of CML $Q/R$. This completes the proof of Lemma 9.
\smallskip\\
\textbf{Theorem 1}.  \textit{The following statements are
equivalent for an arbitrary non-associ\-ati\-ve CML $Q$:}

\textit{1) $Q$ satisfies the maximum condition for subloops;}

\textit{2) if  $Q$ contains a centrally nilpotent subloop of class
$n$, then all its subloops of this type are finitely generated;}

\textit{3) if  $Q$ contains a centrally solvable subloop of class
$s$, then all its subloops of this type are finitely generated;}

\textit{4) at least one maximal associative subloop of  $Q$ is
finitely generated;}

\textit{5) non-invariant associative subloops of  $Q$ are finitely
generated;}

\textit{6) invariant subloops of  $Q$ are finitely generated.}
\smallskip\\
\textbf{Proof.} Let us suppose that CML $Q$ is non-periodic. It
follows from Lemma 4 that subloop $Q^3$ belongs to the center of
CML $Q$. If $H$ is a centrally nilpotent subloop of class $n$
either a centrally solvable subloop of class $s$,  or a maximal
associative subloop, or a non-invariant associative subloop, or an
invariant subloop, then subloop $<H, Q^3>$ will be of this type
too. Therefore it follows from the justice of one of the
statements 2) - 6) of the theorem that abelian group $Q^3$ is
finitely generated. Then it decomposes into a direct product $Q^3
= R \times S$, where $R$ is an abelian group without torsion, $S$
is a finite abelian group [7]. It is obvious that CML $Q/R$ is
periodic. Then by Lemma 5 it is locally finite.

If CML $Q$ satisfies one of the conditions 2) - 6) of  theorem,
then by Lemma 9 CML $Q/R$ satisfies this condition as well. Then
all centrally nilpotent subloops of class $n$ either all centrally
solvable subloops of class $s$,  or at least one maximal
associative subloop, or all non-invariant associative subloops, or
all invariant subloops are respectively finite in CML $Q/R$.
Therefore by Lemma 6 CML $Q/R$ satisfies the minimum condition for
subloops in any case. The center of CML $Q/R$ is finite. Then by
2) of Lemma 6 CML $Q/R$ is finite. Therefore CML $Q$ is finitely
generated and by Lemma 1, the condition 1) holds  in it. It proves
the implications $2) \rightarrow 1), 3) \rightarrow 1), 4)
\rightarrow 1), 5) \rightarrow 1), 6) \rightarrow 1)$. The case
when CML $Q$ is periodic is contained in the proof of previous
case. As the implications $1) \rightarrow 2),  1) \rightarrow 3),
1) \rightarrow 4), 1) \rightarrow 5), 1) \rightarrow 6)$ are
obvious, the theorem is proved.
\smallskip\\
\textbf{Theorem 2}. \textit{The following statements are
equivalent for an arbitrary non-associ\-ati\-ve CML $Q$ with the
multiplication group $\frak M$:}

\textit{1) $Q$ satisfies the maximum condition for subloops;}

\textit{2)  $\frak M$ is finitely generated;}

\textit{3)  $\frak M$ satisfies the maximum condition for
subgroups;}

\textit{4) all invariant subgroups of $\frak M$ is finitely
generated;}

\textit{5) at least one maximal abelian subgroup of $\frak M$ is
finitely generated;}

\textit{6) if  $\frak M$ contains a nilpotent subgroup of class
$n$, then all its subgroups of this type are finitely generated;}

\textit{7) if  $\frak M$ contains a solvable subgroup of class
$s$, then all its subgroups of this type are finitely generated.}
\smallskip\\
\textbf{Proof.} If CML $Q$ satisfies the condition 1), then it is
finitely generated, and by [3] the associator loop $Q'$ is finite.
By Lemma 2 the inner mapping group $I(Q)$ of  $Q$ is also finite.
It is show in [8] that the relation

$$\frak M(G/G') \cong \frak M(G)/<I(G), \textbf
M(G')>,\eqno{(2)}$$ holds  in an arbitrary CML $G$, where $\textbf
M(G')$ denotes a subgroup of group $\frak M(G)$, generated by the
set $\{L(a) \vert a \in G'\}$. It is obvious that group $<I(Q),
\textbf M(Q')>$ is finitely generated in our case. As the abelian
group $\frak M(Q/Q')$ is finitely generated, then it follows from
(2) that group $\frak M$ is finitely generated as well.
Consequently, $1) \rightarrow 2)$.

If the group $\frak M$ is finitely generated, then by Lemma 3 it
is nilpotent. It is known (for instance, see [7]) that the maximum
condition for subgroups holds  in such groups.

Let $Z(Q)$ be the center of an arbitrary CML $Q$, $\{Z(\frak M)\}$
be the upper central series of its multiplication group $\frak
M(Q)$. Then

$$Z(Q) \cong Z(\frak M). \eqno{(3)}$$ Indeed, if $\varphi \in
Z(\frak M)$, then $\varphi L(x) = L(x)\varphi$ for any $x \in Q$.
Further, $\varphi L(x)y = L(x)\varphi y, \varphi (xy) = x\varphi
y$. Let $y = 1$. Then $\varphi x = x\varphi 1, \varphi x =
L(\varphi 1)x, \varphi = L(\varphi 1)$. Now, using the equality
$\varphi(xy) = x\varphi y$ we obtain that $xy\cdot \varphi 1 =
\varphi (xy) = x\cdot\varphi y = x\cdot\varphi (y\cdot 1) = x\cdot
y\varphi 1$. Consequently, if $\varphi \in Z(\frak M)$, then
$\varphi = L(a)$ and $a \in Z(Q)$. Conversely, let $a \in Z(Q)$.
Then $a\cdot xy = ax\cdot y, L(a)L(y)x = L(y)L(a)x, L(a)L(y) =
L(y)L(a)$. It follows from the definition of group $\frak M$ that
$L(a) \in Z(\frak M)$. Finally, if $a, b \in Z(Q)$, then the
homomorphism (3) follows from the equalities $a\cdot bx = ab\cdot
x, L(a)L(b)x = L(ab)x, L(a)L(b) = L(ab)$.

In order to prove the implication $3) \rightarrow 1)$ we use the
relation

$$\frak M/Z_2(\frak M) \cong \frak M(Q/Z(Q)),\eqno{(4)}$$ taking
place in an arbitrary CML [3]. By Lemma 2 the group $\frak
M/Z(\frak M)$ is periodic. Then the group $\frak M/Z_2(\frak M)$,
as an homomorphic image of group $\frak M/Z(\frak M)$, is also
periodic. If the group $\frak M$ satisfies the maximum condition
for subgroups, then the center $Z(\frak M)$ and by (3), also the
center $Z(Q)$, are finitely generated. By Lemma 3 the group $\frak
M$ is nilpotent. Then the group $\frak M/Z_2(\frak M)$ is also
nilpotent and, as it is periodic, then is finite. Hence it follows
from (4) that CML $Q/Z(Q)$ is also finite. Therefore CML $Q$ is
finitely generated and by Lemma 1, the condition 1) holds in it.
Consequently, $3) \rightarrow 1)$.

Let us now suppose that the group $\frak M$ is non-periodic. By
Lemma 2 the group $\frak M/Z(\frak M)$ is locally finite. It
$\alpha$ is an element of infinite order in $\frak M$, then
$\alpha^n \in Z(\frak M)$ for a certain natural number $n$. We
denote by $\frak R$ the subgroup of group $\frak M$, generated by
all elements of form $\alpha^n$. It is obvious that the abelian
group $Z(\frak M)$ is finitely generated if the group $\frak M$
satisfies one of the conditions 4) - 7). Then $Z(\frak M) = \frak
N \times \frak S$, where $\frak N$ is a finitely generated abelian
group without torsion, $\frak S$ is a finite abelian group [7] and
$\frak R = \frak N$. As $\frak N \cap \frak S = \{1\}$, then
$Z(\frak M)/\frak N = (\frak N \times \frak S) \cong \frak S$. By
Lemma 2 the group $\frak M/Z(\frak M)$ is locally finite. It
follows from the relation $\frak M/Z(\frak M) \cong (\frak M/\frak
N)/(Z(\frak M)/\frak N)$ that group $\frak M/\frak N$ is the
extension of the finite group $Z(\frak M)/\frak N$ by locally
finite group $\frak M/Z(\frak M)$. Therefore the group $\frak
M/\frak N$ is locally finite.

By Lemma 2 the commutator group $\frak M'$ is locally finite and
as group $\frak N$ is without torsion, then

$$\frak N \cap \frak M' = \{1\}.$$ Let  either condition 4), or 5), or
6), or 7) hold  in group $\frak M$. By analogy with the proof of
Lemma 9 we can show that in group $\frak M/\frak N$ a condition
analogue with  either conditions 4), or 5), or 6), or 7) holds. We
have already shown that group $\frak M/\frak N$ is locally finite.
Then either all invariant subgroups, or at least one maximal
abelian subgroup, or all nilpotent of class $n$ subgroups, or all
solvable of class $s$ subgroups are finite respectively in $\frak
M/\frak N$. The group $Z(\frak M)/\frak N$ is finite, then it
follows from the relation

$$\frak M/Z(\frak M) \cong (\frak M/\frak N)/(Z(\frak M)/\frak
N)$$ that in group $\frak M/Z(\frak M)$ there holds the same
condition as in group $\frak M/\frak N$. Further, it follows from
the relation

$$\frak M/Z_2(\frak M) \cong (\frak M/Z_1(\frak M))/(Z_2(\frak
M)/Z_1(\frak M)) = (\frak M/Z_1(\frak M))/Z(\frak M/Z_1(\frak
M))$$ that  $\frak M/Z_2(\frak M)$ satisfies the same condition as
group $\frak M/Z_1(\frak M)$, and it follows from (4) that $\frak
M(Q/Z(Q))$ satisfies this condition as well, i.e. either all its
invariant subgroups are finite, or at least one maximal abelian
subgroup is finite, or all its nilpotent subgroups of class $n$
are finite, or all its solvable subgroups of class $s$  are
finite. In such a case, by Lemma 7 the group $\frak M(Q/Z(Q))$
satisfies the minimum condition for subgroups. It is obvious that
the center of group $\frak M(Q/Z(Q))$ is finite. Then by  2) of
Lemma 7 the group $\frak M(Q)/Z(Q)$ is finite, and consequently,
the CML $Q/Z(Q)$ is also finite. The center $Z(\frak M)$ of group
$\frak M$ is finitely generated, then it follows from (3) that the
center $Z(Q)$ of  $Q$ is finitely generated, too. Then the CML $Q$
is finitely generated and by Lemma 1 it satisfies the condition
1). Consequently, if the group $\frak M$ is non-periodic the
implications $4) \rightarrow 1), 5) \rightarrow 1), 6) \rightarrow
1), 7) \rightarrow 1)$ hold.

The case when group $\frak M$ is periodic is proved by analogy for
$\frak N = \{1\}$. Further, as the implications $3)\rightarrow 4),
3)\rightarrow 5), 3)\rightarrow 6), 3)\rightarrow 7)$ are obvious,
the theorem is proved.
\smallskip\\
\textbf{Theorem 3}. \textit{The following conditions are
equivalent for an arbitrary non-associ\-ati\-ve commutative
Moufang $ZA$-loop $Q$ with the multiplication group $\frak M$:}

\textit{1)  $Q$ satisfies the maximum condition for subloops;}

\textit{2) if  $Q$ contains a non-invariant (respect. invariant)
centrally nilpotent subloop of class $n$, then at least one
maximal non-invariant (respect. invariant) centrally nilpotent
subloop of class $n$  is finitely generated;}

\textit{3) if  $Q$ contains a non-invariant (respect. invariant)
centrally solvable subloop of class $s$, then at least one maximal
non-invariant (respect. invariant) centrally solvable subloop of
class $s$  is finitely generated;}

\textit{4) if  $Q$ contains a non-invariant (respect. invariant)
centrally nilpotent subloop of class $n$, then it satisfies the
maximum condition for non-invariant (respect. invariant) centrally
nilpotent subloops of class $n$;}

\textit{5) if  $Q$ contains a non-invariant (respect. invariant)
centrally solvable subloop of class $s$, then it satisfies the
maximum condition for non-invariant (respect. invariant) centrally
solvable subloops of class $s$;}

\textit{6) center $Z(Q)$ of CML $Q$ is finitely generated;}

\textit{7) group $\frak M$ is finitely generated;}

\textit{8) if  $\frak M$ contains a non-invariant (respect.
invariant) nilpotent subgroup of class $n$, then at least one
maximal non-invariant (respect. invariant) nilpotent subgroup of
class $n$  is finitely generated;}

\textit{9) if  $\frak M$ contains a non-invariant (respect.
invariant) solvable subgroup of class $s$, then at least one
maximal non-invariant (respect. invariant) solvable subgroup of
class $s$  is finitely generated;}

\textit{10) if  $\frak M$ contains a non-invariant (respect.
invariant)  nilpotent subgroup of class $n$, then it satisfies the
maximum condition for non-invariant (respect. invariant) nilpotent
subgroups of class $n$;}

\textit{11) if  $\frak M$ contains a non-invariant (respect.
invariant) solvable subgroup of class $s$, then it satisfies the
maximum condition for non-invariant (respect. invariant) solvable
subgroups of class $s$;}

\textit{12) center $Z(\frak M)$ of group $\frak M$ is finitely
generated.}
\smallskip\\
\textbf{Proof.} The implication $1) \rightarrow 2), 1) \rightarrow
3), 1) \rightarrow 4), 1) \rightarrow 5), 1) \rightarrow 6)$ are
obvious. If $H$ is a non-invariant (respect. invariant) centrally
nilpotent of class $n$ (or centrally solvable of class $s$)
subloop of CML $Q$, then subloop $<N, Z(Q)>$ will be of this type
too. Therefore by Lemma 1 the implications $2) \rightarrow 6), 3)
\rightarrow 6), 4) \rightarrow 6), 5) \rightarrow 6)$ hold. Let us
now suppose that the condition 6) holds  in CML $Q$ and let $R$ be
an invariant subloop, defined in Lemma 9. By  4) of Lemma 9 the
center $Z(Q/R)$ is finitely generated and periodic. If follows
from Lemma 5 that $Z(Q/R)$ is finite, and it follows from Lemma 8
that the CML $Q/R$ satisfies the minimum condition for subloops.
As the center $Z(Q/R)$ is finite, then by  2) of Lemma 6, CML
$Q/R$ is finite. By its construction, subloop $L$ is finitely
generated, therefore CML $Q$ is also finitely generated and the
justice of condition 1) follows from Lemma 1. Consequently, the
conditions 1), 2), 3), 4), 5), 6) are equivalent.

The equivalence of conditions 7), 8), 9), 10), 11), 12) is proved
by analogy, using  1), 2) of theorem 2. Finally, the equivalence
of conditions 6), 12) follows from (3). Therefore the theorem is
fully proved.

\smallskip
Tiraspol State University

str. Iablochkin 5

Chi\c{s}in\u{a}u, MD-2069

Moldova

E-mail: aliona2010@yahoo.md; sandumn@yahoo.com
\end{document}